\documentclass[10pt]{amsart}

\usepackage{graphicx}
\usepackage{amsmath}
\usepackage{amssymb}

\numberwithin{equation}{section}
\def\RR{{\mathbb R}}
\def\MM{{\mathbb M}}

\def\wto{\rightharpoonup}

\def\eps{\varepsilon}

\def\ell{l}
\def\mwto{\stackrel{*}{\wto}}
\def\bwto{\stackrel{b}{\wto}}

\def\mint{{{\bf-}\!\!\!\!\!\!\hspace{-.1em}\int}}
\def\limsup{\mathop{\overline{\lim}}}
\def\liminf{\mathop{\underline{\lim}}}
\def\mwto{\stackrel{*}{\wto}}
\def\supp{{\rm supp}}
\def\essinf{\mathop{\rm ess\  inf}}
\def\exp{{\rm e}}

\def\endproof{$\blacksquare$}

\newtheorem{theorem}{Theorem}[section]
\newtheorem{lemma}[theorem]{Lemma}
\newtheorem{proposition}[theorem]{Proposition}
\newtheorem{corollary}[theorem]{Corollary}

\theoremstyle{definition}
\newtheorem{definition}[theorem]{Definition}

\theoremstyle{remark}
\newtheorem{remark}[theorem]{Remark}

\title{Lower semicontinuity via $W^{1,q}$-quasiconvexity}

\author{Jean-Philippe Mandallena}
\address{Laboratoire MIPA (Math\'ematiques, Informatique, Physique et Applications)\newline UNIVERSITE DE NIMES, Site des Carmes, Place Gabriel P\'eri, 30021 N\^\i mes, France.}
\email{jean-philippe.mandallena@unimes.fr}

\keywords{Weak lower semicontinuity, $W^{1,q}$-quasiconvexity, Young measures, equi-integrability, localization principle}

\subjclass[2010]{49J45 (49J10, 74G65)}

\begin{document}

\begin{abstract}
We isolate a general condition, that we call ``localization principle", on the integrand $L:\MM\to[0,\infty]$, assumed to be continuous, under which $W^{1,q}$-quasiconvexity with $q\in[1,\infty]$ is a sufficient condition for $I(u)=\int_\Omega L(\nabla u(x))dx$ to be sequentially weakly lower semicontinuous on $W^{1,p}(\Omega;\RR^m)$ with $p\in]1,\infty[$. Some applications are given.
\end{abstract}

\maketitle


\section{Introduction}

\subsection{The main result} Let $m,N\geq 1$ be two integers, let $\Omega\subset\RR^N$ be a bounded open set with Lipschitz boundary, let $\MM:=\MM^{m\times N}$, where $\MM^{m\times N}$ denotes the space of all real $m\times N$ matrices. Let $p\in]1,\infty[$, let $L:\MM\to[0,\infty]$ be a continuous function and let $I:W^{1,p}(\Omega;\RR^m)\to[0,\infty]$ be defined by
$$
I(u):=\int_\Omega L(\nabla u(x))dx.
$$
In \cite{ball-murat84} Ball and Murat  introduced the concept of $W^{1,q}$-quasiconvexity for $q\in[1,\infty]$, i.e., $L$ is $W^{1,q}$-quasiconvex if and only if
$$
\int_YL\left(\nabla u(y)\right)dy\geq L(\xi)\hbox{ for all }u\in l_\xi+W^{1,q}_0(Y;\RR^m)
$$
with $l_\xi(y):=\xi y$ and $Y:=]-{1\over 2},{1\over 2}[^N$,  and proved (see \cite[Corollary 3.2]{ball-murat84}) that $W^{1,p}$-quasiconvexity is a necessary condition for $I$ to be sequentially weakly lower semicontinuous (swlsc) on $W^{1,p}(\Omega;\RR^m)$, i.e., when 
$$
u_n\wto u\hbox{ in }W^{1,p}(\Omega;\RR^m)\hbox{ implies }\liminf_{n\to\infty}I(u_n)\geq I(u).
$$ 
However, proving that $W^{1,p}$-quasiconvexity, or some variant of it, is also sufficient is still an open problem. In this paper we isolate a general condition on $L$ (see  (C$_{p,q}$) in Theorem \ref{Main-Theorem}) under which $W^{1,q}$-quasiconvexity is a sufficient condition for $I$ to be swlsc on $W^{1,p}(\Omega;\RR^m)$. More precisely, our main result is the following.

\begin{theorem}\label{Main-Theorem}
Given $p\in]1,\infty[$ and $q\in[1,\infty]$, assume that $L$ is $W^{1,q}$-quasiconvex and satisfies
\begin{itemize}
\item[(C$_{p,q}$)] for every $\xi\in\MM$ and every $\{v_n\}_n\subset W^{1,p}(Y;\RR^m)$ such that
$$
\left\{
\begin{array}{l}
v_n\wto l_\xi\hbox{ in }W^{1,p}(Y;\RR^m);\\
\displaystyle\sup_n\int_YL(\nabla v_n(y))dy<\infty,
\end{array}
\right.
$$
there exist a subsequence $\{v_n\}_n$ (not relabeled) and $\{w_n\}_n\subset l_\xi+W^{1,q}_0(Y;\RR^m)$ such that
$$
\left\{
\begin{array}{l}
\displaystyle\left|\nabla v_n-\nabla w_n\right|\to0\hbox{ in measure};\\
\displaystyle\left\{L(\nabla w_n)\right\}_n\hbox{ is equi-integrable.}
\end{array}
\right.
$$
\end{itemize}
Then, $I$ is swlsc on $W^{1,p}(\Omega;\RR^m)$.
\end{theorem}

\subsection{Some applications} First of all, as a direct consequence of Theorem \ref{Main-Theorem}, we have

\begin{corollary}\label{Main-Theorem-Bis}
Given $p\in]1,\infty[$, if {\rm(C$_{p,p}$)} holds then $W^{1,p}$-quasiconvexity is a necessary and sufficient condition for $I$ to be swlsc on $W^{1,p}(\Omega;\RR^m)$. 
\end{corollary}

In fact, Acerbi and Fusco (see \cite{acerbi-fusco84}) showed that $W^{1,\infty}$-quasiconvexity is sufficient for $I$ to be swlsc on $W^{1,p}(\Omega;\RR^m)$ provided that $L$ has $p$-growth, i.e., $L(\cdot)\leq \alpha(1+|\cdot|^p)$ for some $\alpha>0$. We remark that the key argument in their proof is in fact the following result, which we call ``localization principle":
\begin{itemize}
\item[(A$_p$)]{for every $\xi\in\MM$ and every $\{v_n\}_n\subset W^{1,p}(Y;\RR^m)$ such that
$$
v_n\wto l_\xi\hbox{ in }W^{1,p}(Y;\RR^m),
$$ 
there exist a subsequence $\{v_n\}_n$ (not relabeled) and $\{w_n\}_n\subset l_\xi+C^\infty_{\rm c}(Y;\RR^m)$ such that
$$
\left\{
\begin{array}{l}
|\nabla v_n-\nabla w_n|\to0\hbox{ in measure }\\
\{|\nabla w_n|^p\}_n\hbox{ is equi-integrable.}
\end{array}
\right.
$$ }
\end{itemize}
Note that (A) is a particular case of the decomposition lemma (for more details see Kristensen \cite{kristensen94} and also Fonseca, M\"uller and Pedregal \cite{fonseca-muller-pedregal98}). Using this ``localization principle" Kinderlehrer and Pedregal (see \cite{K-P92} and also \cite{sychev99}) proved Acerbi-Fusco's theorem by using Young measure theory. Kinderlehrer-Pedregal's approach was extended by Sychev (see \cite{sychev05a}) to the case where $L$ has fast-growth, i.e., $\beta G(|\cdot|)\leq L(\cdot)\leq \alpha(1+G(|\cdot|))$ for some $\alpha,\beta>0$ and some convex function  $G:[0,\infty[\to[0,\infty[$ such that $\lim_{t\to\infty}{tG^\prime(t)/G(t)}=\infty$ and $tG^\prime(t)/G(t)$ is increasing for large $t$. We also remark that the key argument in its proof is still a ``localization principle", more general than (A$_p$), i.e.,
\begin{itemize}
\item[(B)] {for every $\xi\in\MM$ and every $\{v_n\}_n\subset W^{1,p}(Y;\RR^m)$ such that
$$
\left\{
\begin{array}{l}
v_n\wto l_\xi\hbox{ in }W^{1,p}(Y;\RR^m)\\
\sup\limits_n\int_\Omega G(|\nabla v_n(x)|)dx<\infty,
\end{array}
\right.
$$ 
there exist a subsequence $\{v_n\}_n$ (not relabeled) and $\{w_n\}_n\subset l_\xi+C^\infty_{\rm c}(Y;\RR^m)$ such that
$$
\left\{
\begin{array}{l}
|\nabla v_n-\nabla w_n|\to0\hbox{ in measure }\\
\{G(|\nabla w_n|)\}_n\hbox{ is equi-integrable.}
\end{array}
\right.
$$ }
\end{itemize}

\medskip

It is easily seen that (C$_{p,q}$) generalises (A$_p$) and (B) in a natural way, i.e., 
$$
\left\{
\begin{array}{l}
\hbox{if }L \hbox{ has }p\hbox{-growth then (A$_p$) implies (C$_{p,\infty}$)}\\
\hbox{if }L \hbox{ has fast-growth then (B) implies (C$_{p,\infty}$)},
\end{array}
\right.
$$
which makes that Theorem \ref{Main-Theorem} contains Acerbi-Fusco's theorem and Sychev's theorem in the homogeneous case.

Noticing that the validity of (C$_{p,\infty}$) implies the validity of (C$_{p,q}$) for all $q\in[1,\infty]$, it is obvious that from Theorem \ref{Main-Theorem} we can extend Acerbi-Fusco's theorem and Sychev's theorem, to the case where $L$ is $W^{1,q}$-quasiconvex with $q\in[1,\infty]$, as follows.

\begin{corollary}\label{CoroLL0}
Assume that $L$ is $W^{1,q}$-quasiconvex with $q\in[1,\infty]$.
\begin{itemize} 
\item[(i)] If $L$ has $p$-growth then $I$ is swlsc on $W^{1,p}(\Omega;\RR^m)$.
\item[(ii)] If $L$ has fast-growth then $I$ is swlsc on $W^{1,p}(\Omega;\RR^m)$.
\end{itemize}
\end{corollary}

The following results, i.e., Corollary \ref{CoroLL1}(i)-(ii), are elementary consequences of Theorem \ref{Main-Theorem} in the case where $L$ is $q$-coercive and has $q$-growth with $q$ not necessarily equal to $p$ and $q\not\in\{1,\infty\}$. 

\begin{corollary}\label{CoroLL1}
Assume that $L$ is $W^{1,q}$-quasiconvex and is $q$-coercive and has $q$-growth, i.e., $\beta|\cdot|^q\leq L(\cdot)\leq\alpha (1+|\cdot|^q)$ for some $\alpha,\beta>0$.
\begin{itemize}
\item[(i)] If $p\in]1,N[$ and if $q\in]1,p^*]$ with $p^*:={Np\over N-p}$ then $I$ is swlsc on $W^{1,p}(\Omega;\RR^m)$.
\item[(ii)] If $p\in[N,\infty[$ and if $q\in]1,\infty[$ then $I$ is swlsc on $W^{1,p}(\Omega;\RR^m)$.
\end{itemize}
\end{corollary}

More generally, from Theorem \ref{Main-Theorem} we can establish the following result of Kristensen type (see \cite{kristensen97,kristensen99}).

\begin{corollary}\label{CoroLL2}
If $L$ is finite, $W^{1,q}$-quasiconvex and $q$-coercive and if 
$$
\limsup_{|\xi|\to\infty}{L(\xi)\over|\xi|^q}<\infty,
$$
then {\rm(i)} and {\rm(ii)} of Corollary {\rm\ref{CoroLL1}} are satisfied.
\end{corollary}

\begin{remark}
Similarly to Corollary \ref{CoroLL0}, from the proof of Corollary \ref{CoroLL2} in \S 5.2, it is easily seen that if $\limsup_{|\xi|\to\infty}{L(\xi)/|\xi|^p}<\infty$ and if $L$ is finite and $W^{1,q}$-quasiconvex with $q\in[1,\infty]$, then $I$ is swlsc on $W^{1,p}(\Omega;\RR^m)$.  
\end{remark}
In the following result, inspired by the work of Sychev (see \cite{sychev05a}), we introduce conditions on $L$ (see (D$_1$-D$_3$) in Theorem \ref{Appli-Theorem}) under which (C$_{p,p}$) holds for all $p\in]N,\infty[$.

\begin{theorem}\label{Appli-Theorem}
Assume that $L$ is finite and satisfies{\rm:}
\begin{itemize}
\item[(D$_1$)] there exists $\lambda:]1,\infty[\to]0,1[$ such that $\lambda(R)\to 1$ as $R\to\infty$ and
$$
\lim_{R\to\infty}\sup_{|\xi|\geq R}{L(\lambda(R)\xi)\over L(\xi)}=0;
$$
\item[(D$_2$)] there exists $\alpha_1>0$ such that 
$$
L(t\xi)\leq \alpha_1(1+L(\xi))
$$
for all $\xi\in\MM$ and all $t\in[0,1];$
\item[(D$_3$)] for every $\xi\in\MM$, there exist $\eps>0$, $\alpha_{2,\xi}>0$ and  $c_{\xi}>0$ such that
$$
L(\xi+t(\zeta-\xi)+a)\leq \alpha_{2,\xi}(1+L(\zeta))
$$
for all $t\in[0,1]$, all $\zeta\in\MM$ with $|\zeta-\xi|\geq c_{\xi}$ and all $a\in\MM$ with $|a|\leq \eps$.
\end{itemize}
Then {\rm(C$_{p,p}$)} holds for all $p\in]N,\infty[$.
\end{theorem}

The following result, i.e., Corollary \ref{CoroLL4}, which is a consequence of Theorem \ref{Appli-Theorem} and Corollary \ref{Main-Theorem-Bis}, is not contained in Sychev's theorem, see Remark \ref{Intro-ReMark}. (In fact, Corollary \ref{CoroLL4} is a consequence of a more general result, see Proposition \ref{CoroLL3}, whose statement and proof, which follows from Theorem \ref{Appli-Theorem} and Corollary \ref{Main-Theorem-Bis}, are given in \S 5.3.)

\begin{corollary}\label{CoroLL4}
Let $p\in]N,\infty[$  and let $f:[0,\infty[\to[0,\infty[$ be an increasing concave function with the following properties{\rm:}
\begin{itemize}
\item[(p$_1$)] there exists $\theta\in]0,1[$ such that $f(t s)\leq t^\theta f(s)$ for all $s\in[0,\infty[$ and all $t\in[0,1];$  
\item[(p$_2$)] there exists $r\in]0,1[$ such that $\liminf\limits_{t\to\infty}{f(t)\over t^r}\in]0,\infty];$
\item[(p$_3$)] there exists $\eps>0$ such that $\sup\limits_{0\leq t\leq\eps}f(t)<\infty$.
\end{itemize}
If $L$ has $f(|\cdot|)$-exponential-growth, i.e., $\beta {\rm \exp}^{f(|\cdot|)}\leq L(\cdot)\leq \alpha(1+\exp^{f(|\cdot|)})$ for some $\alpha,\beta>0$, then $I$ is swlsc on $W^{1,p}(\Omega;\RR^m)$ if and only if $L$ is $W^{1,p}$-quasiconvex.
\end{corollary}

\begin{remark}\label{Intro-ReMark}
Corollary \ref{CoroLL4} can be applied in cases which are not covered by Sychev's theorem. Indeed, if $L$ has $f(|\cdot|)$-exponential-growth with $f(t)=t^\nu$ where $\nu\in]0,1[$, then $L$ has $G(|\cdot|)$-growth, i.e., $\beta G(|\cdot|)\leq L(\cdot)\leq\alpha(1+G(|\cdot|))$ for some $\alpha,\beta>0$, with $G(t)=\exp^{t^\nu}$. But $L$ has not fast-growth because such a $G$ is not convex, hence Sychev's theorem cannot be applied. However, as the function $f(t)=t^\nu$ satisfies the hypotheses of Corollary \ref{CoroLL4}, for such a $L$ and for $p\in]N,\infty[$, $I$ is swlsc on $W^{1,p}(\Omega;\RR^m)$ if and only if $L$ is $W^{1,p}$-quasiconvex.
\end{remark}

\subsection*{Plan of the paper} Theorem \ref{Main-Theorem} is proved in Section 3. Its proof uses some classical facts on Young measures that we recall in Section 2. (Note that it seems to be difficult to prove Theorem \ref{Main-Theorem} without using Young measure theory.) Theorem \ref{Appli-Theorem} is proved in Section 4. Finally, Section 5 is devoted to the proofs of corollaries \ref{CoroLL1}, \ref{CoroLL2} and \ref{CoroLL4}.


\section{Some facts on Young measures}

Young measures were introduced by Young in 1937 (see \cite{young37}) with the purpose of finding an extension of the class of Sobolev functions for which one-dimensional nonconvex variational problems become solvable. In the context of the multidimensional calculus of variations, Kinderleherer and Pedregal (see \cite{K-P92,K-P94}) and independently Kristensen (see \cite{kristensen94}) were the first to use Young measures for dealing with lower semicontinuity problems. Relaxation and convergence in energy problems were studied for the first time by Sychev via Young measures following a new approach to Young measures that he introduced in \cite{sychev99}. In this section we only recall the ingredients that we need for proving Theorem \ref{Main-Theorem}. For more details on Young measure theory and its applications to the calculus of variations we refer to \cite{pedregal97,pedregal00,sychev04b}.

Let $\mathcal{P}(\MM)$ be the  set of all probability measures on $\MM$, let $C(\MM)$ be the space of all continuous functions from $\MM$ to $\RR$ and let
$$
C_0(\MM):=\Big\{\Phi\in C(\MM):\lim_{|\xi|\to\infty}\Phi(\xi)=0\Big\}.
$$ 

Here is the definition of a Young measure.

\begin{definition}\label{Def-of-Young-Measures}
A family $(\mu_x)_{x\in\Omega}$ of probability measures on $\MM$, i.e., $\mu_x\in\mathcal{P}(\MM)$ for all $x\in\Omega$, is said to be a Young measure if there exists a sequence $\{\xi_n\}_n$ of measurable functions from $\Omega$ to $\MM$ such that
$$
\Phi(\xi_n)\mwto\langle\Phi;\mu_{(\cdot)}\rangle\hbox{ in }L^\infty(\Omega)\hbox{ for all }\Phi\in C_0(\MM)
$$
with $\langle\Phi;\mu_{(\cdot)}\rangle:=\int_{\MM}\Phi(\zeta)d\mu_{(\cdot)}(\zeta)$. In this case, we say that $\{\xi_n\}_n$ generates $(\mu_x)_{x\in\Omega}$ as a Young measure.
\end{definition}

The following lemma makes clear the link between convergence in measure and Young measures. (The proof follows from the definition.)

\begin{lemma}\label{Lemma1-YM}
let $\{\xi_n\}_n$ and $\{\zeta_n\}_n$ be two sequences of measurable functions from $\Omega$ to $\MM$. If $\{\xi_n\}_n$ generates a Young measure and if $|\xi_n-\zeta_n|\to0$ in measure then $\{\zeta_n\}_n$ generates the same Young measure.
\end{lemma}

The following theorem gives a sufficient condition for proving the existence of Young measures (for a proof see \cite{ball89,sychev04b,fonseca-leoni07}).

\begin{theorem}\label{ExiStenCE-Young-measures}
Let $\theta:\MM\to\RR$ be a continuous function such that 
$
\lim_{|\zeta|\to\infty}\theta(\zeta)=\infty
$
and let $\{\xi_n\}_n$ be a sequence of measurable functions from $\Omega$ to $\MM$ such that
$$
\sup_{n}\int_\Omega \theta(\xi_n(x))dx<\infty.
$$
Then,  $\{\xi_n\}_n$ contains a subsequence generating a Young measure.
\end{theorem}

The following two theorems are important in dealing with integral functionals (for proofs see \cite{balder84,sychev99}).

\begin{theorem}[semicontinuity theorem]\label{S-T-YM}
Let $L:\MM\to[0,\infty]$ be a continuous function and let $\{\xi_n\}_n$ be a sequence of measurable functions from $\Omega$ to $\MM$ such that $\{\xi_n\}_n$ generates $(\mu_x)_{x\in\Omega}$ as a Young measure. Then 
$$
\liminf_{n\to\infty}\int_\Omega L(\xi_n(x))dx\geq\int_\Omega \langle L;\mu_x\rangle dx.
$$
\end{theorem}

\begin{theorem}[continuity theorem]\label{C-T-YM}
Let $L:\MM\to[0,\infty]$ be a continuous function and let $\{\xi_n\}_n$ be a sequence of measurable functions from $\Omega$ to $\MM$ such that $\{\xi_n\}_n$ generates $(\mu_x)_{x\in\Omega}$ as a Young measure. Then 
$$
\lim\limits_{n\to\infty}\int_\Omega L(\xi_n(x))dx=\int_\Omega\langle L;\mu_x\rangle dx<\infty
$$
if and only if $\{L(\xi_n)\}_n$ is equi-integrable.
\end{theorem}


\section{Proof of Theorem \ref{Main-Theorem}}

Let $\{u_n\}_n\subset W^{1,p}(\Omega;\RR^m)$ and let $u\in W^{1,p}(\Omega;\RR^m)$ be such that $u_n\wto u$ in $W^{1,p}(\Omega;\RR^m)$. We have to prove that
\begin{equation}\label{P-MT-Eq0}
\liminf_{n\to\infty}I(u_n)\geq I(u).
\end{equation}
\subsection*{Step 1: localization} Without loss of generality we can assume that:
\begin{eqnarray}
&&\|u_n-u\|_{L^p(\Omega;\RR^m)}\to0;\label{P-MT-Eq1}\\
&&\infty>\liminf_{n\to\infty}I(u_n)=\lim_{n\to\infty}I(u_n) \hbox{ and so }\sup_n\int_\Omega L(\nabla u_n(x))dx<\infty.\label{P-MT-Eq2}
\end{eqnarray}
As $u_n\wto u$ \hbox{ in } $W^{1,p}(\Omega;\RR^m)$ we have 
\begin{equation}\label{Add-Eq1}
\sup_n\int_\Omega|\nabla u_n(x)|^pdx<\infty,
\end{equation} 
and so, by Theorem \ref{ExiStenCE-Young-measures}, there exists a family $(\mu_x)_{x\in\Omega}$ of probability measures on $\MM$ such that (up to a subsequence)
\begin{equation}\label{P-MT-Eq3}
\{\nabla u_n\}_n\hbox{ generates }(\mu_x)_{x\in\Omega}\hbox{ as a Young measure}.
\end{equation}
From Theorem \ref{S-T-YM} it follows that
$$
\liminf_{n\to\infty}I(u_n)\geq \int_{\Omega}\langle L;\mu_x\rangle dx
$$
with (because \eqref{P-MT-Eq2} holds) for a.e. $x_0\in\Omega$,
\begin{equation}\label{P-MT-Eq5}
\langle L;\mu_{x_0}\rangle<\infty.
\end{equation}
Thus, to prove \eqref{P-MT-Eq0} it is sufficient to show that for a.e. $x_0\in\Omega$,
\begin{equation}\label{P-MT-Eq5-goal}
\langle L;\mu_{x_0}\rangle\geq L(\nabla u(x_0)).
\end{equation}

\subsection*{Step 2: blow up} From \eqref{P-MT-Eq2} we deduce that there exist $f\in L^1(\Omega;[0,\infty[)$ and a finite positive Radon measure $\lambda$ on $\Omega$ with $|\supp(\lambda)|=0$ such that (up to a subsequence) $L(\nabla u_n)dx\mwto fdx+\lambda$ in the sense of measures and for a.e. $x_0\in\Omega$,
\begin{equation}\label{P-MT-Eq6}
\lim_{r\to0}\lim_{n\to\infty}\mint_{x_0+rY}L(\nabla u_n(x))dx=f(x_0)<\infty
\end{equation}
with $Y:=]-{1\over 2},{1\over 2}[^N$. By the same argument, from \eqref{Add-Eq1} we see that
\begin{equation}\label{Add-Eq2}
\lim_{r\to0}\lim_{n\to\infty}\mint_{x_0+rY}|\nabla u_n(x)|^pdx<\infty.
\end{equation}
As $u\in W^{1,p}(\Omega;\RR^m)$ it follows that $u$ is a.e. $L^p$-differentiable (see \cite[Theorem 3.4.2 p.129]{ziemer89}), i.e., for a.e. $x_0\in\Omega$,
\begin{equation}\label{P-MT-Eq7}
\lim_{r\to0}{1\over r^{N+p}}\big\|u(x_0+\cdot)-u(x_0)-\nabla u(x_0)y\big\|^p_{L^p({r} Y;\RR^m)}=0.
\end{equation}
From \eqref{P-MT-Eq1} we see that (up to a subsequence) for a.e. $x_0\in\Omega$,
\begin{equation}\label{P-MT-Eq8}
|u_n(x_0)-u(x_0)|^p\to0.
\end{equation}
As $C_0(\MM)$ is separable we can assert that for a.e. $x_0\in\Omega$, $x_0$ is a Lebesgue point of $\langle\Phi;\mu_{(\cdot)}\rangle$ for all $\Phi\in C_0(\MM)$, i.e., 
\begin{equation}\label{P-MT-Eq9}
\lim_{r\to 0}\mint_{x_0+rY}\langle\Phi,\mu_x\rangle dx=\langle\Phi,\mu_{x_0}\rangle\hbox{ for all }\Phi\in C_0(\MM).
\end{equation} 
Fix any $x_0\in\Omega$  such that \eqref{P-MT-Eq5}, \eqref{P-MT-Eq6}, \eqref{P-MT-Eq7}, \eqref{P-MT-Eq8} and \eqref{P-MT-Eq9} hold and fix $r_0>0$ such that $x_0+rY\subset\Omega$ for all $r\in]0,r_0]$. For each $n\geq 1$ and each $r\in]0,r_0]$, let $u_n^r\in W^{1,p}(Y;\RR^m)$ and a family  $(\mu_{y}^r)_{y\in Y}$ of probability measures  on $\MM$ be given by
$$
\left\{
\begin{array}{ll}
u_n^r(y):={1\over r}\left(u_n(x_0+ry)-u_n(x_0)\right)\\
\mu_y^r:=\mu_{x_0+ry}.
\end{array}
\right.
$$
Then \eqref{P-MT-Eq6} (resp. \eqref{Add-Eq2}) can be rewritten as
\begin{equation}\label{P-MT-Eq6-Bis}
\lim_{r\to0}\lim_{n\to\infty}\int_{Y}L(\nabla u_n^r(x))dx<\infty\ \hbox{(resp.} \lim_{r\to0}\lim_{n\to\infty}\int_{Y}|\nabla u_n^r(x)|^pdx<\infty\hbox{)}.
\end{equation}
Taking \eqref{P-MT-Eq3} into account it is easy to see that for every $r\in]0,r_0]$, $\{\nabla u_n^r\}_n$ generates $(\mu_y^r)_{y\in Y}$ as a Young measure, i.e.,
\begin{equation}\label{P-MT-Eq10}
\Phi(\nabla u_n^r)\mwto\langle\Phi,\mu_{(\cdot)}^r\rangle\hbox{ in }L^\infty(Y)\hbox{ as }n\to\infty\hbox{ for all }\Phi\in C_0(\MM),
\end{equation}
and using \eqref{P-MT-Eq9} it is clear that
\begin{equation}\label{P-MT-Eq11}
\langle\Phi;\mu_{(\cdot)}^r\rangle\mwto\langle\Phi;\mu_{x_0}\rangle\hbox{ in }L^\infty(Y)\hbox{ as }r\to0\hbox{ for all }\Phi\in C_0(\MM).
\end{equation} 
On the other hand, we have
\begin{eqnarray*}
\|u_{n}^{r}-l_{\nabla u(x_0)}\|^p_{L^p(Y;\RR^m)}&=&\int_{Y}|u_{n}^{r}(y)-l_{\nabla u(x_0)}(y)|^pdy\\
&=&{1\over r^{N+p}}\|u_n(x_0+\cdot)-u_n(x_0)-l_{\nabla u(x_0)}\|^p_{L^p(rY;\RR^m)},
\end{eqnarray*}
and consequently
\begin{eqnarray*}
\|u_{n}^r-l_{\nabla u(x_0)}\|^p_{L^p(Y;\RR^m)}&\leq&{c\over r^{N+p}}\|u_n-u\|^p_{L^p(\Omega;\RR^m)}+{c\over r^{N+p}}|u_n(x_0)-u(x_0)|^p\\
&&+{c\over r^{N+p}}\|u(x_0+\cdot)-u(x_0)-l_{\nabla u(x_0)}\big\|^p_{L^p(rY;\RR^m)}
\end{eqnarray*}
with $c>0$ which only depends on $p$. Using \eqref{P-MT-Eq1}, \eqref{P-MT-Eq8} and \eqref{P-MT-Eq7} we deduce that
\begin{equation}\label{P-MT-Eq12}
\lim_{r\to0}\lim_{n\to\infty}\|u_n^r-l_{\nabla u(x_0)}\|_{L^p(Y;\RR^m)}=0.
\end{equation}
According to \eqref{P-MT-Eq12}, \eqref{P-MT-Eq6-Bis} and \eqref{P-MT-Eq10} together with \eqref{P-MT-Eq11}, by diagonalization there exists a mapping $n\to r_n$ decreasing to $0$ such that
$$
\left\{
\begin{array}{l}
\|u_n^{r_n}-l_{\nabla u(x_0)}\|_{L^p(Y;\RR^m)}\to0\\
\lim\limits_{n\to\infty}\int_{Y}|\nabla u_n^{r_n}(y)|^pdy<\infty, \hbox{ and so }\sup\limits_n\int_Y|\nabla u_n^{r_n}(y)|^pdy<\infty\\
\lim\limits_{n\to\infty}\int_{Y}L(\nabla u_n^{r_n}(y))dy<\infty, \hbox{ and so }\sup\limits_n\int_YL(\nabla u_n^{r_n}(y))dy<\infty\\
\{\nabla u_n^{r_n}\}_n\hbox{ generates }\mu_{x_0}\hbox{ as a Young measure},
\end{array}
\right.
$$
and consequently we have:
\begin{eqnarray}
&&\left\{
\begin{array}{l}
v_n\wto l_{\nabla u(x_0)}\hbox{ in }W^{1,p}(Y;\RR^m)\\
\sup\limits_n\int_YL(\nabla v_n(y))dy<\infty;
\end{array}
\right.\label{Hyp-Cpq}\\
&&\{\nabla v_n\}_n\hbox{ generates }\mu_{x_0}\hbox{ as a Young measure}.\label{Hyp-CT-YM}
\end{eqnarray}
where $v_n:=u_n^{r_n}$. 

\subsection*{Step 3: using (C\boldmath$_{p,q}$\unboldmath)\  and\ \boldmath$W^{1,q}$\unboldmath-quasiconvexity} According to \eqref{Hyp-Cpq}, by (C$_{p,q}$) there exists $\{w_n\}_n\subset l_{\nabla u(x_0)}+W^{1,q}_0(Y;\RR^m)$ such that
$$
\left\{
\begin{array}{l}
|\nabla v_n-\nabla w_n|\to0\hbox{ in measure}\\
L(\nabla w_n)\hbox{ is equi-integrable,}
\end{array}
\right.
$$
hence, by \eqref{Hyp-CT-YM} and  Lemma \ref{Lemma1-YM}, $\{\nabla w_n\}_n$ generates $\mu_{x_0}$ as a Young measure, and, taking \eqref{P-MT-Eq5} into account, from Theorem \ref{C-T-YM} we deduce that
\begin{equation}\label{P-MT-Eq13}
\lim_{n\to\infty}\int_YL(\nabla w_n(y))dy=\langle L;\mu_{x_0}\rangle.
\end{equation}
 As $L$ is $W^{1,q}$-quasiconvex, we have
$$
\int_YL(\nabla w_n(y))dy\geq L(\nabla u(x_0))\hbox{ for all }n\geq 1,
$$
and \eqref{P-MT-Eq5-goal} follows by letting $n\to\infty$ and using \eqref{P-MT-Eq13}. 
$\blacksquare$

\begin{remark}
In case $q=\infty$ the condition of $W^{1,q}$-quasiconvexity is the classical condition of quasiconvexity by Morrey (see \cite{morrey52}).
\end{remark}

\begin{remark}
In fact, we have also proved that if $\{u_n\}_n\subset W^{1,p}(\Omega;\RR^m)$ is such that $\sup_{n}\int_\Omega L(\nabla u_n(x))dx<\infty$ and if $\{\nabla u_n\}_n$ generates $(\mu_x)_{x\in\Omega}$ as a Young measure, then for a.e. $x\in\Omega$, $\mu_x$ is a homogeneous gradient $L$-Young measure centered at $\nabla u(x)$, with $u\in W^{1,p}(\Omega;\RR^m)$, provided that $u_n\wto u$  in  $W^{1,p}(\Omega;\RR^m)$ and (C$_{p,q}$) holds with $q\in[1,\infty]$. Homogeneous gradient $L$-Young measures were introduced and completely characterized by Sychev in \cite{sychev00} where we refer the reader for more details.
\end{remark}

\begin{remark} From the proof of Theorem \ref{Main-Theorem} we can extract the following lower semicontinuity theorem with the biting weak convergence.
\begin{theorem}
Given $p\in]1,\infty[$ and $q\in[1,\infty]$, assume that $L$ is $W^{1,q}$-quasiconvex and satisfies {\rm(C$_{p,q}$)}. Then, for each $\{u_n\}_n\subset W^{1,p}(\Omega;\RR^m)$ and each $u\in W^{1,p}(\Omega;\RR^m)$ such that $u_n\wto u$ in $W^{1,p}(\Omega;\RR^m)$ and $\sup_n\int_\Omega L(\nabla u_n(x))dx<\infty$,  there exists a subsequence $\{u_n\}_n$ (not relabeled) and a family $(\mu_x)_{x\in\Omega}$ of probability measures on $\MM$ such that{\rm:}
\begin{itemize}
\item[(i)] $\{\nabla u_n\}_n$ generates $(\mu_x)_{x\in\Omega}$ as a Young measure{\rm;}
\item[(ii)] $L(\nabla u_n)\bwto\langle L;\mu_(\cdot)\rangle$, where ``$\bwto$" denotes the biting weak convergence{\rm;}
\item[(iii)] $\langle L;\mu_x\rangle\geq L(\nabla u(x))$ for a.a. $x\in\Omega$.
\end{itemize}
\end{theorem}
\end{remark}

For a deeper discussion of weak lower semicontinuity in the sense of biting lemma, see Ball and Zhang \cite{ball-zhang90} (see also \cite[Lemma 3.2]{sychev05a} for a simple proof of the biting lemma).


\section{Proof of Theorem \ref{Appli-Theorem}}

Let $p\in]N,\infty[$, let $\xi\in\MM$ and let $\{v_n\}_n\subset W^{1,p}(Y;\RR^m)$ be such that: 
\begin{eqnarray}
&& v_n\wto l_\xi\hbox{ in }W^{1,p}(Y;\RR^m),\hbox{ and so }\sup_n\int_Y|\nabla v_n(y)|^pdy<\infty;\label{Eq-Appli-1}\\
&&\sup_n\int_YL(\nabla v_n(y))dy<\infty\label{Eq-Appli-3}.
\end{eqnarray}
As $p>N$, \eqref{Eq-Appli-1} implies that, up to a subsequence, 
\begin{equation}\label{Eq-Appli-4}
\|v_n- l_\xi\|_{L^\infty(Y;\RR^m)}\to0.
\end{equation}
\subsection*{Step 1: using the biting Lemma} First, recall Sychev's version of the biting lemma (see \cite[Lemma 3.2]{sychev05a}).
\begin{lemma}\label{biting-lemma}
Let $\{f_n\}_n\subset L^1(Y;[0,\infty[)$ be such that $\sup_n\int_Yf_n(y)dy<\infty$. Then, there exist a subsequence $\{f_n\}_n$ (not relabeled) and $\{M_n\}_n\subset]0,\infty[$ with $M_n\to\infty$ such that $\{f_n\chi_{Y_n}\}_n$ is equi-integrable with $\chi_{Y_n}$ denoting the characteristic function of $Y_n:=\{y\in Y:f_n(y)\leq M_n\}$.
\end{lemma}

Taking \eqref{Eq-Appli-3} into account, from Lemma \ref{biting-lemma}, that we apply with $f_n=L(\nabla v_n)$, we can assert that, up to a subsequence, 
\begin{equation}\label{Eq-Appli-5}
\{L(\nabla v_n)\chi_{Y_n}\}_n\hbox{ is equi-integrable.} 
\end{equation}
Let $\{R_n\}_n$ be given by
$
R_n:=\essinf_{y\in Y\setminus Y_n}|\nabla v_n(y)|.
$
As $L$ is finite and continuous, $Y\setminus Y_n=\{y\in Y:L(\nabla v_n(y))>M_n\}$ and $M_n\to\infty$, we have
$
R_n\to\infty.
$
Let $\{u_n\}_n\subset W^{1,p}(Y;\RR^m)$ be defined by
$$
u_n:=\lambda_nv_n\hbox{ with }\lambda_n:=\lambda(R_n),
$$
where $\lambda:]1,\infty[\to]0,1[$, with $\lambda(R_n)\to1$, is given by (D$_1$). From \eqref{Eq-Appli-1} and \eqref{Eq-Appli-4} we have:
\begin{eqnarray}
&& \|\nabla v_n-\nabla u_n\|_{L^p(Y;\MM^{m\times N})}\to0;\label{Eq-Appli-6bis}\\
&&\|u_n-l_\xi\|_{L^\infty(Y;\RR^m)}\to0.\label{Eq-Appli-6}
\end{eqnarray}
On the other hand, given any $n\geq 1$, $L(\nabla u_n)=L(\lambda_n \nabla v_n)\chi_{Y_n}+L(\lambda_n \nabla v_n)\chi_{Y\setminus Y_n}$, and so $L(\nabla u_n)\leq\alpha_1(1+L(\nabla v_n)\chi_{Y_n})+L(\lambda_n \nabla v_n)\chi_{Y\setminus Y_n}$ by using (D$_2$). But $|\nabla v_n|\geq R_n$ on $Y\setminus Y_n$, hence 
$$
L(\nabla u_n)\leq\alpha_1\left(1+L(\nabla v_n)\chi_{Y_n}\right)+\sup_{|\zeta|\geq R_n}{L(\lambda_n\zeta)\over L(\zeta)}L(\nabla v_n).
$$
Taking \eqref{Eq-Appli-3} and \eqref{Eq-Appli-5} into account and noticing that $\sup_{|\zeta|\geq R_n}{L(\lambda_n\zeta)\over L(\zeta)}\to0$ by (D$_1$) we conclude that 
\begin{equation}\label{Eq-Appli-7}
\{L(\nabla u_n)\}_n\hbox{ is equi-integrable}.
\end{equation}

\subsection*{Step 2: cut-off method} Set
$$
\theta_\xi:=\sup_{|\zeta|\leq |\xi|+c_\xi+\eps}L(\zeta)\hbox{ with }c_{\xi}>0\hbox{ and }\eps>0\hbox{ given by (D$_3$)}.
$$
(Such a $\theta_\xi$ exists because $L$ is finite and continuous.) Fix any $n\geq 1$. Let $\phi_n\in C_c^\infty(Y;[0,1])$ be a cut-off function between $Q_n:=]{\eps_n-1\over2},{1-\eps_n\over 2}[^N$ and $Y$ such that $\|\phi_n\|_{L^\infty(Y)}\leq{2\over\eps_n}$ with 
$$
\eps_n:={\|u_n-l_\xi\|^ {1\over 2}_{L^\infty(Y;\RR^m)}}.
$$
(Note that, by \eqref{Eq-Appli-6}, $\eps_n\to 0$.) Define $w_n\in l_\xi+W^{1,p}_0(Y;\RR^m)$ by
$$
w_n:=l_\xi+\phi_n(u_n-l_\xi).
$$
Then
$$
\nabla w_n=\left\{
\begin{array}{ll}
\nabla u_n&\hbox{on }Q_n\\
\xi+\phi_n(\nabla u_n-\xi)+\nabla\phi_n\otimes(u_n-l_\xi)&\hbox{on }Y\setminus Q_n.
\end{array}
\right.
$$
Setting $C_n:=\{y\in Y\setminus Q_n:|\nabla u_n(y)-\xi|<c_\xi\}$ we have
\begin{equation}\label{Eq-Appli-8}
L(\nabla w_n)\leq L(\nabla u_n)+L(\nabla w_n)\chi_{C_n}+L(\nabla w_n)\chi_{Y\setminus(Q_n\cup C_n)}.
\end{equation}
But $|\xi+\phi_n(\nabla u_n-\xi)+\nabla\phi_n\otimes(u_n-l_\xi)|\leq|\xi|+c_\xi+2\eps_n$ on $C_n$ and $\eps_n\to 0$, hence $|\xi+\phi_n(\nabla u_n-\xi)+\nabla\phi_n\otimes(u_n-l_\xi)|\leq|\xi|+c_\xi+\eps$ on $C_n$, and so
\begin{equation}\label{Eq-Appli-9}
L(\nabla w_n)\chi_{C_n}\leq\theta_\xi\chi_{Y\setminus Q_n}.
\end{equation}
Moreover, as $0\leq\phi_n\leq 1$, $|\nabla u_n-\xi|\geq c_\xi$ on $Y\setminus(Q_n\cup C_n)$ and $|\nabla\phi_n\otimes(u_n-l_\xi)|\leq 2\eps_n\to0$, from (D$_3$) we see that
\begin{equation}\label{Eq-Appli-10}
L(\nabla w_n)\chi_{Y\setminus(Q_n\cup C_n)}\leq\alpha_{2,\xi}(1+L(\nabla u_n)).
\end{equation}
Combining \eqref{Eq-Appli-9} and \eqref{Eq-Appli-10} with \eqref{Eq-Appli-8} we deduce that for every $n\geq 1$,
$$
L(\nabla w_n)\leq \theta_\xi\chi_{Y\setminus Q_n}+(\alpha_{2,\xi}+1)(1+L(\nabla u_n)).
$$
Taking \eqref{Eq-Appli-7} into account and noticing that $|Y\setminus Q_n|\to0$, we deduce that 
$$
\{L(\nabla w_n)\}_n\hbox{ is equi-integrable.}
$$
On the other hand, it is easy to see that there exists $K>0$, which only depends on $p$, such that
$$
|\nabla w_n-\nabla u_n|^p\leq K\left(|\xi|^p+|\nabla v_n|^p\right)\chi_{Y\setminus Q_n}+K\eps_n^p.
$$
Taking \eqref{Eq-Appli-1} into account and recalling that $|Y\setminus Q_n|\to 0$ and $\eps_n\to0$, we deduce that $\|\nabla w_n-\nabla u_n\|_{L^p(Y;\MM^{m\times N})}\to0$, and so $\|\nabla v_n-\nabla w_n\|_{L^p(Y;\MM^{m\times N})}\to0$ by combining with \eqref{Eq-Appli-6bis}. It follows that
$$
|\nabla v_n-\nabla w_n|\to0\hbox{ in measure,}
$$
and the proof is complete. $\blacksquare$

\section{Proof of Corollaries \ref{CoroLL1}, \ref{CoroLL2} and \ref{CoroLL4}}

\subsection{Proof of Corollary \ref{CoroLL1}} In each case, it is sufficient to prove that (C$_{p,q}$) holds. Let $\xi\in\MM$ and let $\{v_n\}_n\subset W^{1,p}(Y;\RR^m)$ be such that:
\begin{eqnarray}
&& v_n\wto l_\xi\hbox{ in }W^{1,p}(Y;\RR^m);\label{HYPOTH1}\\
&& \sup_n\int_YL(\nabla v_n(y))dy<\infty.\label{HYPOTH2}
\end{eqnarray}
From \eqref{HYPOTH1} we deduce that
\begin{eqnarray}
&&\sup_n\|v_n\|_{W^{1,p}(Y;\RR^m)}<\infty\label{HYPOTH3}
\end{eqnarray}
and, combining \eqref{HYPOTH2} with the fact that $L$ is $q$-coercive, we obtain
\begin{eqnarray}
&&\sup_n\|\nabla v_n\|_{L^q(Y;\MM)}<\infty.\label{HYPOTH4}
\end{eqnarray}
On the other hand, using Sobolev imbeddings, we have:
\begin{itemize}
\item[(i)] if $p\in]1,N[$ and $q\in]1,p^*]$ then there exists $C_1>0$ such that 
\begin{equation}\label{HYPOTH5}
\|v_n\|_{L^q(Y;\RR^m)}\leq C_1\|v_n\|_{W^{1,p}(Y;\RR^m)}\hbox{ for all }n\geq 1;
\end{equation}
\item[(ii)] if $p\in[N,\infty[$ and $q\in]1,\infty[$ then there exists $C_2>0$ such that 
\begin{equation}\label{HYPOTH6}
\|v_n\|_{L^q(Y;\RR^m)}\leq C_2\|v_n\|_{W^{1,p}(Y;\RR^m)}\hbox{ for all }n\geq 1.
\end{equation}
\end{itemize}
Thus, combining \eqref{HYPOTH3} and \eqref{HYPOTH4} with \eqref{HYPOTH5} or \eqref{HYPOTH6}, we see that, in each case, $\sup_n\|v_n\|_{W^{1,q}(Y;\RR^m)}<\infty$, and so (up to a subsequence) $v_n\wto l_\xi$ in $W^{1,q}(Y;\RR^m)$ by considering \eqref{HYPOTH1}. From (A$_q$) (which corresponds to (A$_p$), stated in \S 1.1, with $p=q$) we deduce that there exist a subsequence $\{v_n\}_n$ (not relabeled) and $\{w_n\}_n\subset l_\xi+W^{1,q}_0(Y;\RR^m)$ such that $|\nabla v_n-\nabla w_n|\to 0$ in measure and $\{|\nabla w_n|^q\}_n$ is equi-integrable, and so $\{L(\nabla w_n)\}_n$ is equi-integrable because $L$ has $q$-growth, which completes the proof. \endproof

\subsection{Proof of Corollary \ref{CoroLL2}} In each case, we have to prove that (C$_{p,q}$) is satisfied. Arguing as in the proof of Corollary \ref{CoroLL2}, given $\xi\in\MM$ and $\{v_n\}_n\subset W^{1,p}(Y;\RR^m)$ verifying \eqref{HYPOTH1} and \eqref{HYPOTH2}, there exist a subsequence $\{v_n\}_n$ (not relabeled) and $\{w_n\}_n\subset l_\xi+W^{1,q}_0(Y;\RR^m)$ such that $|\nabla v_n-\nabla w_n|\to 0$ in measure and 
\begin{eqnarray}
\{|\nabla w_n|^q\}_n\hbox{ is equi-integrable.}\label{HYPOTH7}
\end{eqnarray}
We are thus reduced to prove that
\begin{eqnarray}
\{L(\nabla w_n)\}_n\hbox{ is equi-integrable.}\label{HYPOTH8}
\end{eqnarray}
As $\theta:=\limsup_{|\xi|\to\infty}{L(\xi)\over|\xi|^q}<\infty$ there exists $R>0$ such that
\begin{eqnarray}
L(\xi)\leq (\theta+1)|\xi|^q\hbox{ for all }|\xi|> R.\label{HYPOTH9}
\end{eqnarray}
As $L$ is finite and continuous we have
\begin{eqnarray}
M:=\sup_{|\xi|\leq R}L(\xi)<\infty.\label{HYPOTH10}
\end{eqnarray}
Set $Y_n:=\{y\in Y:|\nabla w_n(y)|\leq R\}$. Then, given any $n\geq 1$, $L(\nabla w_n)=L(\nabla w_n)\chi_{Y_n}+L(\nabla w_n)\chi_{Y\setminus Y_n}$, where $\chi_E$ denotes the characteristic function of the set $E$. Using \eqref{HYPOTH9} and \eqref{HYPOTH10} we see that $L(\nabla w_n)\leq M\chi_{Y_n}+(\theta+1)|\nabla w_n|^q\chi_{Y\setminus Y_n}$, hence
$$
L(\nabla w_n)\leq M+(\theta+1)|\nabla w_n|^q,
$$
and \eqref{HYPOTH8} follows from \eqref{HYPOTH7}. \endproof

\subsection{Proof of Corollary \ref{CoroLL4}} It is a consequence of the following proposition whose proof, which follows from Theorem \ref{Appli-Theorem} and Corollary \ref{Main-Theorem-Bis}, is given below. 

\begin{proposition}\label{CoroLL3}
Fix $p\in]N,\infty[$ and assume that $L$ has $(F,\gamma)$-exponential-growth, i.e., $\beta {\rm \exp}^{\gamma F(\cdot)}\leq L(\cdot)\leq \alpha(1+\exp^{F(\cdot)})$ for some $\alpha,\beta>0$, some $\gamma\geq 1$ and some Borel measurable function $F:\MM\to[0,\infty[$, with $F$ satisfying the following conditions{\rm:}
\begin{itemize}
\item[(d$_1$)] there exists $\theta\in]0,\infty[$ such that $F(t\xi)\leq t^\theta F(\xi)$ for all $\xi\in\MM$ and all $t\in[0,1];$ 
\item[(d$_2$)] there exists $r\in]0,\infty[$ such that $\liminf\limits_{|\xi|\to\infty}{F(\xi)\over|\xi|^r}\in]0,\infty];$ 
\item[(d$_3$)] there exists $\eps>0$ such that $\sup\limits_{|\xi|\leq\eps}F(\xi)<\infty;$ 
\item[(d$_4$)] $F(\zeta-\xi)\leq F(\zeta)+F(\xi)$ for all $\zeta,\xi\in\MM;$
\item[(d$_5$)] $F(\zeta+\xi)\leq \sqrt{\gamma}\big(F(\zeta)+F(\xi)\big)$ for all $\zeta,\xi\in\MM$.
\end{itemize}
Then, $I$ is swlsc on $W^{1,p}(\Omega;\RR^m)$ if and only if $L$ is $W^{1,p}$-quasiconvex.
\end{proposition}

According to Proposition \ref{CoroLL3}, it is sufficient to prove that $F:\MM\to[0,\infty[$ defined by $F(\xi):=f(|\xi|)$ satisfies (d$_1$-d$_4$) and (d$_5$) with $\gamma=1$. It is obvious that (p$_i$) implies (d$_i$) for $i=1,2,3$. Thus, we only need to show that (d$_4$) and (d$_5$) holds with $\gamma=1$. As $f:[0,\infty[\to[0,\infty[$ is concave we can assert that $f$ is subadditive, i.e., 
\begin{equation}\label{Subadditive-Property}
f(s+t)\leq f(s)+f(t)\hbox{ for }s,t\in[0,\infty[.
\end{equation}
Fix $\xi,\zeta\in\MM$. Then, $|\zeta-\xi|\leq|\zeta|+|\xi|$, hence $f(|\zeta-\xi|)\leq  f(|\zeta|+|\xi|)$ because $f$ is increasing, and so $f(|\zeta-\xi|)\leq f(|\zeta|)+f(|\xi|)$, i.e., $F(\zeta-\xi)\leq F(\zeta)+F(\xi)$ , by using \eqref{Subadditive-Property}. On the other hand, $|\zeta+\xi|\leq|\zeta|+|\xi|$, hence $f(|\zeta+\xi|)\leq  f(|\zeta|+|\xi|)$ because $f$ is increasing, and so $f(|\zeta+\xi|)\leq f(|\zeta|)+f(|\xi|)$, i.e., $F(\zeta+\xi)\leq F(\zeta)+F(\xi)$ , by using \eqref{Subadditive-Property}, and the proof is complete. \endproof

\subsection*{Proof of Proposition \ref{CoroLL3}} It is sufficient to prove that (C$_{p,p}$) holds. To do this, we are going to establish (D$_1$), (D$_2$) and (D$_3$) of Theorem \ref{Appli-Theorem}.

First of all, as $L$ has $(F,\gamma)$-exponential-growth, given any $\lambda:]1,\infty[\to]0,1[$, by using (d$_1$) we have 
$$
{L(\lambda(R)\xi)\over L(\xi)}\leq{\alpha\left(1+\exp^{F(\lambda(R)\xi)}\right)\over\beta\exp^{\gamma F(\xi)}}\leq{\alpha\over\beta}\left(\exp^{-\gamma F(\xi)}+\exp^{\left(\lambda^\theta(R)-\gamma\right) F(\xi)}\right)
$$
for all $R>1$ and all $\xi\in\MM$ with $\alpha,\beta>0$, $\gamma\geq 1$ and $\theta\in]0,\infty[$ given by (d$_1$). From (d$_2$), setting either $\delta:={1\over 2}\liminf_{|\xi|\to\infty}{F(\xi)\over|\xi|^r}$ if $\liminf_{|\xi|\to\infty}{F(\xi)\over|\xi|^r}\in]0,\infty[$ or $\delta=1$ if $\liminf_{|\xi|\to\infty}{F(\xi)\over|\xi|^r}=\infty$ (and so $\lim_{|\xi|\to\infty}{F(\xi)\over|\xi|^r}=\infty$), with $r\in]0,\infty[$,  
we can assert that there exists $R_\delta>1$ such that 
\begin{equation}\label{HypOtHeSis-d2}
F(\xi)\geq \delta|\xi|^r\hbox{ for all }|\xi|\geq R_\delta.
\end{equation}
For each $R\geq R_\delta$, as $\lambda^\theta(R)-\gamma\leq \lambda^\theta(R)-1<0$, by using \eqref{HypOtHeSis-d2}, we see that
$$
{L(\lambda(R)\xi)\over L(\xi)}\leq{\alpha\over\beta}\left(\exp^{-\gamma\delta R^r}+\exp^{\delta\left(\lambda^\theta(R)-1\right) R^r}\right)
$$ 
whenever $|\xi|\geq R$, and  consequently 
$$
\sup_{|\xi|\geq R}{L(\lambda(R)\xi)\over L(\xi)}\leq{\alpha\over\beta}\left(\exp^{-\gamma\delta R^r}+\exp^{\delta\left(\lambda^\theta(R)-1\right) R^r}\right)
$$
for all $R\geq R_\delta$, which shows that (D$_1$) holds with $\lambda(R)=\big(1-R^{-{r\over 2}}\big)^{1\over\theta}$.

By the fact that $L$ has $(F,\gamma)$-exponential-growth, for every $\xi\in\MM$ and every $t\in[0,1]$,  by using (d$_1$) we have
$$
L(t\xi)\leq \alpha(1+\exp^{F(\xi)})\leq \alpha(1+\exp^{\gamma F(\xi)})\leq\alpha\max\left\{1,{1\over\beta}\right\}(1+L(\xi)),
$$
which shows that (D$_2$) holds with $\alpha_1=\alpha\max\big\{1,{1\over\beta}\big\}$.

Finally, as $L$ has $(F,\gamma)$-exponential-growth, by using (d$_1$), (d$_4$) and (d$_5$) we see that
$$
L(\xi+t(\zeta-\xi)+a)\leq\alpha\left(1+\exp^{\sqrt{\gamma} F(a)}\exp^{2\gamma F(\xi)}\exp^{\gamma F(\zeta)}\right)
$$
for all $\xi,\zeta, a\in\MM$ and all $t\in[0,1]$. From (d$_3$) there exists $\eps>0$ such that
$$
M:=\sup_{|a|\leq \eps}F(a)<\infty,
$$
and consequently, for every $\xi\in\MM$, we have
$$
L(\xi+t(\zeta-\xi)+a)\leq\alpha\max\left\{1,{\exp^{\sqrt{\gamma} M}\exp^{2\gamma F(\xi)}\over\beta}\right\}\left(1+L(\zeta)\right)
$$
for all $t\in[0,1]$, all $\zeta\in\MM$ and all $a\in\MM$ with $|a|\leq \eps$, which shows that (D$_3$) holds with $\alpha_{2,\xi}=\alpha\max\big\{1,{\exp^{\sqrt{\gamma} M}\exp^{2\gamma F(\xi)}\over\beta}\big\}$. \endproof

\medskip

\subsection*{Acknowledgments} I gratefully acknowledges the many comments of M. A. Sychev during the preparation of this paper, and for the lectures on ``Young measures and weak convergence theory" that he gave at the University of N\^imes during may-june 2011.

The author also wishes to thank the ``R\'egion Languedoc Roussillon" for financial support through its program ``d'accueil de personnalit\'es \'etrang\`eres" which allowed to welcome M. A. Sychev, from the Sobolev Institute for Mathematics in Russia, in the Laboratory MIPA of the University of N\^imes.

\bibliographystyle{alpha}

\end{document}